\documentclass{amsart}
\usepackage{amsmath, amssymb}
\newtheorem{Theorem}{Theorem}[section]
\newtheorem{Cor}[Theorem]{Corollary}
\newtheorem{Lemma}[Theorem]{Lemma}
\newtheorem{Proposition}[Theorem]{Proposition}
\theoremstyle{definition}
\newtheorem{Definition}[Theorem]{Definition}
\theoremstyle{remark}
\newtheorem{rem}[Theorem]{Remark}
\numberwithin{equation}{section}
\newcommand{\R}{\mathbb R}
\newcommand{\N}{\mathbb N}

\newcommand{\C}{\mathcal{C}}
\newcommand{\F}{\mathcal{F}}
\newcommand{\p}{\mathcal{P}}
\newcommand{\A}{\mathcal{A}}
\newcommand{\D}{\mathcal{D}}
\newcommand{\id}{Id_{]0;1[}}
\begin{document}

\title{Deformed algebras: examples and application to Lax equations}
\author{Jean-Pierre Magnot}
%\author{}%

\address{Lyc\'ee Jeanne d'Arc - Avenue de Grand Bratagne - F-63000 Clermont-Ferrand}
\email{jean-pierr.magnot@ac-clermont.fr}
%\email{}
%\email{}%

%\commby{}%
% ----------------------------------------------------------------

\begin{abstract}
We develop here a concept of deformed algebras through three examples and an application. Deformed algebras are obtained from a fixed algebra by deformation along a family of indexes, through formal series. We show how the example of deformed algebra used in \cite{Ma2013} is only an example among others, and how they often give rise to regular Fr\"olicher Lie groups. Then, we show how such deformed algebras arise in a formal integration of Lax equations through time scaling. The infinite dimensional groups under consideration enables to state the uniqueness of the formal solutions, their smooth dependance under perturbation, and to study some of the symmetries.      
\end{abstract}
\maketitle

\tableofcontents
\section*{Introduction}
In \cite{Ma2013}, an algebra and a group of formal series of operators is described in order to rewrite the integration of the KP hierarchy in a non formal way. One of the main advances of this work is to get a (non formal) principal bundle where the concept of holonomy makes sense rigorously. The geometric objects under consideration are diffeological or Fr\"olicher groups, which are regular in the sense that the exponential map exists and is smooth. 

Diffeological spaces, first described in the 80's by Souriau and his coworkers \cite{Don, Igdiff, Les, Sou} are generalizations of manifolds that enables differential geometry without charts. Independently, Fr\"olicher spaces give a more rigid framework, that also generalize the notion of manifolds \cite{CN, FK, KM}. The comparison of the two frameworks has been made independently in \cite{Ma2006} and in \cite{Wa}, see e.g. \cite{Ma2013}. 

The aim of this paper is to enlarge the spectrum of application of such formal series of algebras. First by summarizing and enlarging the theoretical results of \cite{Ma2013}, secondly giving three examples:

- one that enlarges straightway the example of \cite{Ma2013}: in $q-$deformation, indexes are monomials $q^n$, and in the example that we develop the base algebra is $A=Cl(M,E),$ the algebra of classical pseudo-differential operators. It is not an enlargeable Lie algebra, but the formal series $$G = 1+\sum_{n >0} q^n A$$ is a regular Fr\"olicher Lie group (section \ref{opd}). 

- one can also consider indexes that are manifolds: this is the case of path-like or cobordism-like deformations (section \ref{pc}) or tensor algebras (section \ref{tens}) . 

Finally, we give an application of  $q-$deformation in section \ref{Laxequation}: A Lax equation \cite{Lax} is a formally integrable equation of the type : 
$$\partial_t L = [P,L]$$
where $P,L$ are in most cases differential, pseudo-differential, or difference operators. This equation integrates heuristically  as a classical equation on a group of matrices: there should have an unique solution, up to the initial value $L(0)$, given by $$L(t) = Ad_{Exp P (t)} L(0).$$
Unfortunately, very often, the operator   $ Exp P (t)$ exists only at a formal level. We propose to apply a scaling $$t \mapsto qt$$ to the time variable. The operator $P(t)$ is changed into an operator $P_q(t)$ which is a monomial of order $1$ in the $q-$variable. This allows the machinery of $q-$deformed operators: the algebras considered are now Lie algebras of (smooth) regular Lie groups. As a simple consequence, we get smoothness of the unique solution $L_q(t)$ with respect to $P(t)$ and $L(0);$ another consequence is that a class of symmetries of the $q-$deformed Lax equation obey also a Lax-type equation, and hence can be described with the same techniques, again in a rigorous sense in the class of $q-$deformed solutions and depending smoothly of the parameters.
A natural open question is the limit $q \rightarrow 1$, where the solutions of the initial equation $\partial_t L = [P,L]$ should be recovered. A work in progress studies this limit for classical examples. 
  
\section*{Acknowledgements} I would like to thank Professor Ambar Sengupta for stimulating discussions on the topics of cobordism, that influenced the corresponding section of this paper.   
	
\section{Regular Fr\"olicher Lie groups of series of unbounded operators} 

\subsection{Diffeological spaces and
Fr\"olicher spaces}
\label{1.1}

\begin{Definition} Let $X$ be a set.

\noindent $\bullet$ A \textbf{plot} of dimension $p$ (or $p$-plot)
on $X$ is a map from an open subset $O$ of $\R^{p}$ to $X$.

\noindent $\bullet$ A \textbf{diffeology} on $X$ is a set $\p$
of plots on $X$ such that, for all $p\in\N$,

- any constant map $\R^{p}\rightarrow X$ is in $\p$;

- Let $I$ be an arbitrary set; let $\{f_{i}:O_{i}\rightarrow X\}_{i\in I}$
be a family of maps that extend to a map $f:\bigcup_{i\in I}O_{i}\rightarrow X$.
If $\{f_{i}:O_{i}\rightarrow X\}_{i\in I}\subset\p$, then $f\in\p$.

- (chain rule) Let $f\in\p$, defined on $O\subset\R^{p}$. Let $q\in\N$,
$O'$ an open subset of $\R^{q}$ and $g$ a smooth map (in the usual
sense) from $O'$ to $O$. Then, $f\circ g\in\p$.

\vskip 6pt $\bullet$ If $\p$ is a diffeology $X$, $(X,\p)$ is
called \textbf{diffeological space}.

\noindent Let $(X,\p)$ et $(X',\p')$ be two diffeological spaces,
a map $f:X\rightarrow X'$ is \textbf{differentiable} (=smooth) if
and only if $f\circ\p\subset\p'$. \end{Definition}

We now introduce Fr\"olicher spaces.

\begin{Definition} $\bullet$ A \textbf{Fr\"olicher} space is a triple
$(X,\F,\C)$ such that

- $\C$ is a set of paths $\R\rightarrow X$,

- A function $f:X\rightarrow\R$ is in $\F$ if and only if for any
$c\in\C$, $f\circ c\in C^{\infty}(\R,\R)$;

- A path $c:\R\rightarrow X$ is in $\C$ (i.e. is a \textbf{contour})
if and only if for any $f\in\F$, $f\circ c\in C^{\infty}(\R,\R)$.

\vskip 5pt $\bullet$ Let $(X,\F,\C)$ et $(X',\F',\C')$ be two
Fr\"olicher spaces, a map $f:X\rightarrow X'$ is \textbf{differentiable}
(=smooth) if and only if $\F'\circ f\circ\C\in C^{\infty}(\R,\R)$.
\end{Definition}

Any family of maps $\F_{g}$ from $X$ to $\R$ generate a Fr\"olicher
structure $(X,\F,\C)$, setting \cite{KM}:

- $\C=\{c:\R\rightarrow X\hbox{ such that }\F_{g}\circ c\subset C^{\infty}(\R,\R)\}$

- $\F=\{f:X\rightarrow\R\hbox{ such that }f\circ\C\subset C^{\infty}(\R,\R)\}.$

One easily see that $\F_{g}\subset\F$. This notion will be useful
in the sequel to describe in a simple way a Fr\"olicher structure.

A Fr\"olicher space, as a differential space, carries a natural topology,
which is the pull-back topology of $\R$ via $\F$. In the case of
a finite dimensional differentiable manifold, the underlying topology
of the Fr\"olicher structure is the same as the manifold topology. In
the infinite dimensional case, these two topologies differ very often.

\vskip 12pt

One can remark, if $X$ is a Fr\"olicher space, we define a natural
diffeology on $X$ by \cite{Ma2006}, see e.g. \cite{Ma2013, Wa}: 
$$
\p(\F)=
\coprod_{p\in\N}\{\, f\hbox{ p-
paramatrization on } X; \, \F \circ f \in C^\infty(O,\R) \quad \hbox{(in
the usual sense)}\}.$$
With this construction, we also get a natural diffeology when
$X$ is a Fr\"olicher space. In this case, one can easily show the following:
\begin{Proposition}\cite{Ma2006}, see e.g. \cite{Ma2013, Wa}
Let$(X,\F,\C)$
and $(X',\F',\C')$ be two Fr\"olicher spaces. A map $f:X\rightarrow X'$
is smooth in the sense of Fr\"olicher if and only if it is smooth for
the underlying diffeologies. \end{Proposition}

Thus, we can also state:

\begin{tabular}{ccccc}
smooth manifold  & $\Rightarrow$  & Fr\"olicher space  & $\Rightarrow$  & Diffeological space\tabularnewline
\end{tabular}

\subsection{Regular Fr\"olicher groups}

Let $(G, \F, \C)$ be a Fr\"'olicher space which is a group  such that
the group law and the inversion map are smooth. These laws are also smooth
for the underlying diffeology. Then, following \cite{Les}, this is
possible as in the case of manifolds to define a tangent space and
a Lie algebra $\mathfrak{g}$ of $G$ using germs of smooth maps.
Let us precise the algebraic, diffeological and Fr\"olicher structures of $\mathfrak{g}.$

\begin{Proposition}
Let $\mathfrak{g} = \{ \partial_t c(0) ; c \in \C \hbox{ and } c(0)=e_G \}$
be the space of germs of paths at $e_G.$
\begin{itemize}
	\item Let $(X,Y) \in \mathfrak{g}^2,$ $X+Y = \partial_t(c.d)(0)$  where $c,d \in \C ^2,$ $c(0) = d(0) =e_G ,$ 
	$X = \partial_t c(0)$ and $Y = \partial_t d(0).$ 
	\item Let $(X,g) \in \mathfrak{g}\times G,$ $Ad_g(X) = \partial_t(g c g^{-1})(0)$  where $c \in \C ,$ $c(0) =e_G ,$ 
	and $X = \partial_t c(0).$  
	\item Let $(X,Y) \in \mathfrak{g}\times G,$ $[X,Y] = \partial_t( Ad_{c(t)}Y)$   where $c \in \C ,$ $c(0) =e_G ,$ 
	$X = \partial_t c(0).$
\end{itemize}
All these operations are smooth and thus well-defined. 
\end{Proposition} 

The basic properties remain globally the same as in the case of Lie groups, and the prrofs are similar replacing charts by plots of the underlying diffeologies. (see e.g. \cite{Les} for further details) 

\begin{Definition} A Fr\"olicher group $G$ with Lie algebra $\mathfrak{g}$
is called \textbf{regular} if and only if there is a smooth map \[
Exp:C^{\infty}([0;1],\mathfrak{g})\rightarrow C^{\infty}([0,1],G)\]
 such that $g(t)=Exp(v(t))$ if and only if $g$ is the unique solution
of the differential equation \[
\left\{ \begin{array}{l}
g(0)=e\\
\frac{dg(t)}{dt}g(t)^{-1}=v(t)\end{array}\right.\]
 We define \begin{eqnarray*}
exp:\mathfrak{g} & \rightarrow & G\\
v & \mapsto & exp(v)=g(1)\end{eqnarray*}
 where $g$ is the image by $Exp$ of the constant path $v.$ \end{Definition}

The classical setting for infinite dimensional differential geometry requires the model topological vector space
to be complete or Mac-Key complete. One of the reasons for this choice 
is to ensure the existence of the integral of a path over a 
compact interval. This means that the choice of an adaquate topology is 
necessary. For vector spaces, the basis for such a study can be found 
in \cite{KM}, when the properties of the so-called ``convenient vector spaces''
are given. We have to remark that a vector space for which addition and scalar multiplication
are compatible with a given Fr\"olicher structure needs only 
a topological structure to become a convenient vector space.
In order to circumvent these topological considerations, and adapting the     
terminology of regular Lie groups to vector spaces (which are viewed as abelian Lie groups), we set:

\begin{Definition}
Let $(V,\F, \C)$ be a Fr\"olicher vector space, i.e. a vector space $V$ equipped with a Fr\"'olicher structure compatible
with the vector space addition and the scalar multiplication. $(V,\F, \C)$ is \textbf{regular} if there is a smooth map 
$$ \int_0^{(.)} : C^\infty([0;1];V) \rightarrow C^\infty([0;1],V)$$ such that $\int_0^{(.)}v = u$ if and only if $u$ is the unique solution of 
the differential equation
\[
\left\{ \begin{array}{l}
u(0)=0\\
u'(t)=v(t)\end{array}\right. .\]

\end{Definition}

This definition is of course fulfilled if $V$ is a complete locally convex topological vector space, equipped with its natural Fr\"'olicher structure.

\begin{Definition}
Let $G$ be a Fr\"olicher Lie group with Lie algebra $\mathfrak{g}.$ Then, $G$ is \textbf{regular with regular Lie algebra}
if both $G$ and $\mathfrak{g}$ are regular.
\end{Definition}

The first known example is the following \cite{Ma2013}:
\begin{Proposition} \label{omo} Let $(G_{n})_{n\in\N}$ be a sequence of Banach
Lie groups, increasing for $\supset,$ and such that the inclusions
are Lie group morphisms. Let $G=\bigcap_{n\in\N}G_{n}.$ Then, $G$
is a Fr\"olicher regular Lie group with regular Lie algebra $\mathfrak{g}=\bigcap_{n\in\N}\mathfrak{g}_{n}.$
\end{Proposition} 
Let us notice that there exists non regular Fr\"olicher Lie groups, see the appendix, where as there is no example of Fr\'echet Lie group that has been proved to be non regular \cite{KM}.
We now turn to  key results from \cite{Ma2013}: 

\begin{Theorem} \label{regulardeformation}
Let $(A_n)_{n \in \N^*} $ be a sequence of complete locally convex (Fr\"olicher)
vector spaces which are regular, 
equipped with a graded smooth multiplication operation
on $ \bigoplus_{n \in \N^*} A_n ,$ i.e. a multiplication such that 
$A_n .A_m \subset A_{n+m},$ smooth with respect to the corresponding Fr\"olicher structures.
Let us assume that:

Then, the set 
$$1 + \A = \left\{ 1 + \sum_{n \in \N^*} a_n | \forall n \in \N^* , a_n \in A_n \right\} $$
is a Fr\"olicher Lie group, with regular  Fr\"olicher Lie algebra
$$\A= \left\{ \sum_{n \in \N^*} a_n | \forall n \in \N^* , a_n \in A_n \right\}.$$
Moreover, the exponential map defines a bijection $\A \rightarrow 1+\A.$  
\end{Theorem}
\begin{Theorem}\label{exactsequence}
Let 
$$ 1 \rightarrow K \underrightarrow{i} G \underrightarrow{p}  H \rightarrow 1 $$
be an exact sequence of Fr\"olicher Lie groups, such that there is a smooth section $s : H \rightarrow G,$ and such that 
the trace diffeology from $G$ on $i(K)$ coindides with the push-forward diffeology from $K$ to $i(K).$
We consider also the corresponding sequence of Lie algebras
$$ 0 \rightarrow \mathfrak{k} \underrightarrow{i'} \mathfrak{g} \underrightarrow{p}  \mathfrak{h} \rightarrow 0 . $$
Then, 
\begin{itemize}
\item The Lie algebras $\mathfrak{k}$ and $\mathfrak{h}$ are regular if and only if the
Lie algebra $\mathfrak{g}$ is regular;
\item The Fr\"olicher Lie groups $K$ and $H$ are regular if and only if the Fr\"olicher Lie group $G$ is regular.
\end{itemize}

\end{Theorem}
We mimick and extend the procedure used in \cite{Ma2013} 
for groups of formal pseudo-differential operators applied to KP equations. 
This method is based on theorem 
\ref{regulardeformation}. For this, we need a family of indexes $I$ that are $\N$-graded, 
equipped with a graded 
composition $*$ that is defined on a domain of $I \times I,$ which is associative, and for which there is a finite number of indexes at the same order.  

\begin{Definition}
Let $I$ as above.
Let $A_i$ be a sequence of regular Fr\"olicher vector spaces.
Let $$\mathcal{A} = \left\{ \sum_{i = 1}^{+\infty} a_i | a_i \in \mathcal{A}_i \right\}.$$ The Fr\"olicher vector space $\mathcal{A}$ is called Fr\"olicher $I-$graded regular algebra if and only if it is equipped with a multiplication, associative and distributive with respect addition, smooth for the induced fr\"olicher structure, such that 
$$A_i . A_j \left\{ \begin{array}{ll} \subset A_{i*j} & \hbox{ if } i*j \hbox{ exists} \\
= 0 & \hbox{ otherwise. } \end{array} \right.$$
\end{Definition}

\begin{Proposition} \label{AI}
Let $\mathcal{A}$ be a Fr\"olicher $I-$graded regular algebra. It is the Lie algebra of the Fr\"olicher regular Lie group $1+\mathcal{A}.$
\end{Proposition}

\noindent
\textbf{Proof.}
This is a straight application of theorem \ref{regulardeformation}
\qed 
Each example will be an application of the following theorem:
\begin{Theorem} \label{extension}
Let $\mathcal{A} = \bigoplus_{i \in I} \mathcal{A}_i$ be a Fr\"olicher $I-$graded regular algebra. Let $G$ be a regular Fr\"olicher Lie group, acting on $\mathcal{A}$ componentwise. Then,  $$G \oplus A$$ is a regular Fr\"olicher Lie group.
\end{Theorem}

\noindent
\textbf{Proof.}
Considering the exact sequence 
$$ 0 \rightarrow 1 + \A \rightarrow G \oplus \A \rightarrow G \rightarrow 0$$
there  is a (global) slice $G \rightarrow G \oplus \{0_\A \}$ so that Theorem \ref{exactsequence} applies . \qed
 \subsection{Examples of $q-$deformed pseudo-differential operators} \label{opd}
In our work of Lax-type equations, we use the following group from \cite{Ma2013}:Let $E$ be a smooth vector bundle over a
compact manifold without boundary M. We denote by  $ Cl(M, E) $ (resp.  $ Cl^k (M, E)
$) the space of
 classical pseudo-differential operators (resp.
classical pseudo-differential operators of order k) acting on smooth
sections of $E$. We denote by $Cl^*(M,\mathbb{C}^n)$,
 $Cl^{0,*}(M,\mathbb{C}^n)$ the groups of
the units of the algebras $Cl(M,\mathbb{C}^n)$ and
$Cl^{0}(M,\mathbb{C}^n)$.
Notice that $Cl^{0,*}(M,\mathbb{C}^n)$ is a CBH Lie group, and belong 
to a wider class of such groups that is studied in \cite{Glo}.

.
\begin{Definition}
Let $t$ be a formal parameter. 
We define the algebra of formal series 
$$Cl_t(M,E) = \left\{ \sum_{t \in \N^*} t^k a_k | \forall k \in \N^*, a_k \in Cl^k(M,E) \right\}.$$
\end{Definition}
This is obviously an algebra, graded by the order (the valuation) into the variable  $t.$ Thus, setting
$$ \A_n = \left\{ t^n a_n | a_n \in Cl^n(M,E)\right\} ,$$
we can set $\A = Cl_t(M,E)$ and state the followinc consequence of Theorem \ref{regulardeformation}:

\begin{Cor}
The group $1 + Cl_t(M,E)$ is a regular Fr\"olicher Lie group with regular 
Fr\"olicher Lie algebra $Cl_t(M,E).$
\end{Cor}
Let $Cl^{0,*}(M,E)$ be the Lie group of invertible pseudo-differential operators of order 0. This group is known to be a regular Lie group since Omori, but the most efficient proof is actually in \cite{Glo}, to our knowledge.
We remark a short exact sequence of Fr\"olicher Lie groups:
$$ 0 \rightarrow 1 + Cl_q(M,E) \rightarrow Cl^{0,*}(M,E) + Cl_q(M,E) \rightarrow Cl^{0,*}(M,E) \rightarrow 0,$$  
which satisfies the conditions of Theorem \ref{extension}. Thus, we have the following:
\begin{Theorem}
$Cl^{0,*}(M,E) + Cl_q(M,E)$ is a regular Fr\"olicher Lie group with Lie algebra $Cl^{0}(M,E) + Cl_q(M,E).$
\end{Theorem}

Let us now go deeper in the constructions allowed by Theorem \ref{extension}:

\begin{Theorem} \label{nonclassicpdo}
$$Cl_q^{\infty,*}(M,E) = \left\{ A = \sum_{k = 0}^{+\infty} q^k a_k \in Cl(M,E)
\left[  \left[ q \right] \right] | a_0 \in Cl^{0;*}(M,E) \right\}$$
is a Fr\"olicher Lie group with Lie algebra
$$Cl_q^{\infty}(M,E) = \left\{ A = \sum_{k = 0}^{+\infty} q^k a_k \in PDO(M,E)
\left[  \left[ q \right] \right] | a_0 \in Cl^{0}(M,E) \right\}.$$
 
\end{Theorem}

\begin{rem} One could also develop a similar example, 
which could stand as a generalized version, 
with log-polyhomogeneous pseudo-differential operators or with other algebras of non classical operators. 
These examples are not developed here in order to avoid some too long lists of examples 
constructed in the same spirit.
\end{rem}

\subsection{Path-like and Cobordism-like deformations} \label{pc}
Let us now consider the set 
$$ Gr = \coprod Gr_n$$
where $Gr_m$ is the set of n-dimensional connected oriented manifolds $M$, possibly with boundary, where the boundaries $\partial M$ are separated into two disconnected parts: the initial part $\alpha(M)$ and the final part $\beta(M).$
Then, we have a composition law $*$ by the relation:

\begin{Definition}
Let $m \in \N^*.$ Let $M, M' \in Gr_m.$ Then $M'' = M*M' \in Gr_n$ exists if 
\begin{enumerate}
\item $\alpha(M) = \beta(M') \neq \emptyset,$ up to diffeomorphism
\item $\alpha(M'')= \alpha(M')$
\item $\beta(M'') = \beta(M)$
\item $M''$ cuts into two pieces $M'' = M \cup M'$ with $M \cap M' = \alpha(M) = \beta(M').$
\end{enumerate}
\end{Definition} 

This composition, that we call \textbf{cobordism composition}, extends naturally to embedded manifolds:
\begin{Definition}
Let $N$ be a smooth manifold.
$$ Gr(N) = \coprod_{m \in \N^*} \coprod_{M \in Gr_m} Emb(M,N).$$
where the notation $Emb(M,N)$ denotes the smooth manifold of smooth embeddings of $M$ into $N.$
\end{Definition}

Notice that $Gr(N)$ is naturally a smooth manifold, since $Emb(M,N)$ is a smooth manifold \cite{KM}, and that $*$ is obviously smooth because it is smooth in the sense of the underlying diffeologies. 
When we only consider manifolds without boundary (in this case,  cobordism composition is not defined), these spaces are called non linear grassmanians in the litterature, which explains the notations.
Let us now turn to $q-$deformed groups and algebras:
\begin{Definition}
\begin{itemize}
\item Let $ I =  \left( Gr \times \N^*\right) \coprod (\emptyset, 0),$ graded by the second component.  Assuming $\emptyset$ as a neutral element for $*$, we extend the cobordism composition into a composition, also noted $*$, defined as:
$$ (M,p) * (M',p') = (M*M', p+ p')$$ when $M*M'$ is defined.
\item
Let $ I(N) =  \left( Gr(N) \times \N^*\right) \coprod (\emptyset, 0),$ graded by the second component.  Assuming $\emptyset$ as a neutral element for $*$, we extend the cobordism composition into a composition, also noted $*$, defined as:
$$ (M,p) * (M',p') = (M*M', p+ p')$$ when $M*M'$ is defined.
\item Let $m \in \N^* .$ We note by $I_m$ and $I_m(N)$ the set of indexes based on $Gr_m.$
\end{itemize}
\end{Definition}
  Let $A$ be a regular Fr\"olicher algebra. Let $m \in \N^*.$ Let 
$$ \A_{I_m} = \left\{ \sum_{(M,n) \in I_m} q^na_{M,n} | a_{M;n} \in A \right\} $$
and let 
$$ \A_{I_m}(N) = \left\{ \sum_{(\phi,n) \in I_m(N)} q^na_{\phi,n} | a_{M;n} \in A \right\}. $$
\begin{Theorem} \label{cob}
Let $\Gamma \subset \coprod_{m \in \N^*} I_m$, resp. $\Gamma(N) \subset \coprod_{m \in \N^*} I_m(N)$, be a finitely generated family of indexes, stable under $*$. Then 
$$ \A_\Gamma = \left\{ \sum_{(M,n) \in \Gamma} q^na_{M,n} | a_{M;n} \in A \right\} ,$$ 
resp. $$ \A_\Gamma(N) = \left\{ \sum_{(\phi,n) \in \Gamma(N)} q^na_{\phi,n} | a_{\phi;n} \in A \right\} ,$$
is a regular algebra (in the sense of regular vector spaces). If $(\emptyset, 0) \notin \Gamma,$ $$ 1_A + \A_\Gamma$$
is a Lie group with Lie algebra $A_\Gamma.$
  If $(\emptyset, 0) \in \Gamma,$ for eacg Fr\"olicher Lie group $G $ with Lie algebra $\mathfrak{g}$ such that $G \subset A^*$ smoothly, $$ G \oplus \A_{\Gamma - (\emptyset, 0)}$$
is a regular Fr\"olicher Lie group with Lie algebra $\mathfrak{g} \oplus A_\Gamma.$
Moreover, the results are the same replacing $\Gamma$ by $\Gamma(N).$
\end{Theorem} 	

\noindent
\textbf{Proof.}
Since $\Gamma$ is finitely generated, in the (possibily infinite sum) $\sum_{(M,n) \in \Gamma} q^na_{M,n},$ each power $q^n$ has only a finite number of $A-$coefficients since there is only a finite number of possible indexes for each $q^n.$ So that, Proposition \ref{AI} and Theorem \ref{extension} apply. The same arguments are also valid when replacing $\Gamma$ by $\Gamma(N).$ \qed

\begin{Cor}
$1+ \A_{I_1}$ is a regular Fr\"olicher Lie group with Lie algebra $\A_{I_1}.$  
\end{Cor}  

\noindent
\textbf{Proof.} $Gr_1$ has 5 elements: $S^1$, $[0;1[,$ $]0,1[;$ $]0;1]$ and $[0;1],$ where the open brackets are understood as lack of boundary, and that the initial point is at $0$, the final one is at $1.$ This set is obviously stable under $*$, so that Theorem \ref{cob} applies. \qed    
\subsection{Tensor product algebras} \label{tens}

Let $I$ be a family of manifolds, stable and finitely generated under cartesian product.  
Cartesian product is the composition law. We remark that it is graded with respect to the dimension of the manifold. A standard singleton can be added to $I$ as a neutral element of dimension $0$.
Let $J$ be a family of finite rank vector bundles over the family of manifolds $I,$ stable under tensor product. The scalar field can be understood as a neutral element of dimension $0.$
   By the way, Proposition \ref{AI} applies to the following algebras :
	
	- algebras of smooth sections of the finite rank vector bundles of $J$
	
	- algebras of operators acting on these sections. 
	
	When $I = \left\{(S^1)^n |n \in \N^* \right\},$ $J = \left\{(S^1)^n \times \mathbb{C} |n \in \N^* \right\},$
	and when the algebras under consideration are $Cl(S^1,\mathbb{C})^{\otimes n}$, we recognize a framework in the vincinity of the example given at the end of \cite{Ma2008}.  
\section{Application to integration via time scaling of a Lax equation} \label{Laxequation}
A Lax equation is an equation of the type:
\begin{equation}\label{lax1}
\left\{ \begin{array}{lll} \partial_t L(t) & = & \left[ P(t), L(t) \right] \\
 L(0) & = & \hbox{fixed operator (initial value)} \end{array} \right.
\end{equation}
where $t \mapsto P(t)$ is a given is a smooth path of operators. If the path $P$
is a smooth path of the Lie algebra $\mathfrak{g}$ of a regular Lie group $G,$
if $G$ acts on a Fr\'echet algebra of operators $\mathcal{B}$ that contains $L(0),$ 
the path \begin{equation} \label{solution1} L(t) = Ad_{Exp_GP(t)} L(0) \end{equation}
is a solution of equation \ref{lax1}. Unfortunately, if $P$ is e.g. a path of differential operators or order bigger than 1, 
there is no known inclusion into an adequate Lie algebra $\mathfrak{g}$ (see e.g. the explanations given in \cite{Ma2008}). one of the biggest known classes of Lie algebras that integrate into a Lie group are particular Lie algebras of classical pseudo-differential operators of order 1, that integrate into algebras of Fourier integral operators of order 0.

\subsection{Formal integration via time scaling}
We only assume that both $P$ and $L$ are in a fixed Fr\'echet algebra $\mathcal{A}$ with unit element,
or in a $c^\infty-$algebra if one prefers to work in the convenient setting \cite{KM}. 
Let us now build a corresponding Lax equation in $\mathcal{A}[[q]].$
We consider the paths $P(qt)$ and $L(qt)$ obtained by time scaling $t \mapsto qt. $
Then, \begin{eqnarray*} \partial_t L(qt) & = & q(\partial_tL)(qt) \\
&=&\left[ qP(qt), L(qt) \right] \end{eqnarray*} 
for a fixed parameter $q.$ We note by $L_q(t)=L(qt)$ and by $P_q(t) = qP(qt).$
We get the following equation:
\begin{equation}\label{lax2}
\left\{ \begin{array}{lll} \partial_t L_q(t) & = & \left[ P_q(t), L_q(t) \right] \\
 L(0) & = & \hbox{fixed operator (initial value) in } \mathcal{A} \end{array} \right.
\end{equation}
Let $val_q$ be the valuation of formal series in $\mathcal{A}[[q]]$ with respect to the $q$ variable.
We remark that $val_qL_q=0$ and $val_qP_q=1. $ We note by $\mathcal{A}[[q]]_{>0}$ the ideal made of formal series S such that $val_qS>0.$
\begin{Theorem}
The equation \ref{lax2} has an unique solution in $\mathcal{A}[[q]]$, with explicit solution:
$$ L_q(t)= exp(P_q)(t).L(0).\left(exp(P_q)(t)\right)^{-1}$$
where the the map $exp$ is the group exponential $\mathcal{A}[[q]]_{>0} \rightarrow Id + \mathcal{A}[[q]]_{>0}.$
\end{Theorem}

The proof is a straightforward consequence of basic results on Lie groups.

\subsection{Smooth dependence of the solutions on the parameters}
The serie $exp(P_q)(t),$ read as 
$$  exp(P_q)(t) =  \sum_{i = 0}^\infty a_i(q) $$ where 
$$a_i (q) = \int_{t\geq s_1 \geq ... \geq s_i\geq 0} \left[ \prod_{j = 1}^i P_q(s_j)\right] (ds)^i
$$

Notice that the coefficients $a_i (q)$ are series in the $q$ parameter if $P_q$ is a serie in the $q$ parameter. But in our definition, $P_q$ is a degree 1 monomial in $q.$ So that, each $a_i$ is a monomial of degree $i$ in the variable $q$ which makes the sum consistent in $\mathcal{A}[[q]]$ and  we have the following: 
\begin{Lemma}
The maps $$(t,q,P)  \mapsto a_i (q)$$
are smooth maps $$\R^2 \times \mathcal{A} \rightarrow \mathcal{A}[q]_i$$
\end{Lemma} 

\noindent
\textbf{Proof.}
The Fr\"olicher structure on $\mathcal{A}[[q]]$ is the one that is generated for the smoothness of each coefficient, and the map $(q,P) \mapsto P_q$ is smooth. By composition with the (smooth) exponential map, we get that each $a_i$ is smooth under the variables $(t,q,P)$. \qed 

\begin{Theorem}
The map $$ (t,q,P,L(0)) \in \R^2 \times \mathcal{A}^2 \mapsto L_q(t) \in \mathcal{A}[[q]]$$
is smooth.  
\end{Theorem}

\noindent
\textbf{Proof.}

\subsection{Symmetries of Lax equations and time scaling}
Let us now look for symmetries of a Lax equation. A symmetry is a path  $S$ of linear invertible operators on $\mathcal{A}$ such that, if $L(t)$ is a solution of \ref{lax1}, $S(t).L(t)$ is also a solution. Assuming smoothness, we shall quickly go into more restricted classes of symmetries. 
 \begin{eqnarray*} \partial_t(S.L)(t) & = & (\partial_t S).L(t) + S. (\partial_t L) (t) \\
& = & (\partial_t S).L(t) + S.[P,L](t) \\
& = & (\partial_t S).L(t) + \left( S.[P,L](t) -[P,S.L](t) \right) + \partial_t(S.L)(t)  
\end{eqnarray*}
from which we get \begin{equation} \label{symmetry1}
(\partial_t S).L(t) = \lbrace ad_P,S \rbrace .L(t)\end{equation}
applying the time scaling, we get, with the obvious notations: 
\begin{equation} \label{symmetry2}
(\partial_t S_q ).L_q(t) =  \left[ ad_{P_q}, S_q\right] .L_q(t)\end{equation}
The map $S \rightarrow S_q$ is an homomorphism from the group of symmetries of (\ref{lax1}) to the group of symmetries of (\ref{lax2}), and it appears to us that there should exist symmetries of (\ref{lax2}) that are not induced from symmetries of (\ref{lax1}).  
\begin{itemize}
\item
The map $\phi : S \mapsto (\partial_t S ).L(t) -  \left[ ad_{P}, S\right] .L(t)$ is linear and the symmetries of (\ref{symmetry1}) are the zeros of $phi.$
The same way, the map $\phi_q : S_q \mapsto (\partial_t S_q ).L_q(t) -  \left[ ad_{P_q}, S_q\right] .L_q(t)$ is linear and the symmetries of (\ref{symmetry2}) are the zeros of $phi_q.$ 
Such a problem appears non relevant to the methods of resolution of this paper, and we leave the question of solving these two equations open. Let us now turn to a special class of solutions. 
\item
Let us now simplify this equation, avoiding the $L_q -$ term. 
Then, we get another Lax-type equation
\begin{equation} \label{symmetry3}
\partial_t S_q  =  \left[ ad_{P_q}, S_q\right] \end{equation}

and we can remark that the operator $ad_{P_q}$ is an inner derivation of $A,$ which is of order 1 in $q$ since $P_q$ is of order 1. Let $In(A)$ be the Lie algebra of inner derivations of $A.$ Let $In_q(A)$ be the $q-$graded algebra of operators spanned by $  qIn_A  ,$ endowed with the push-forward Fr\"olicher structure from $A.$ We get:

\begin{Theorem}
\begin{enumerate}
\item \label{no1}$In_q$ is a smooth regular algebra
\item \label{no2} $Id_A + In_q(A)$ is a regular Fr\"olicher Lie group with Lie algebra $In_q(A).$ \end{enumerate}
\end{Theorem} 

\vskip 12pt
\noindent
\textbf{Proof.}
Let us remark that (\ref{no2}) is a straightforward consequence of (\ref{no1}) and Theorem \ref{regulardeformation}.
Now, we recall that smoothness in $In(A)$ is induced by smoothness in $A.$ Moreover, the inclusion $In(A) \rightarrow C^\infty(A,A)$ is smooth in the Fr\"olicher sense \cite{KM}. So that $In_q(A)$ is a smooth algebra,
where the composition is smooth and bilinear. Finally, by the standard properties of the integral, we get that some paths $ad_{a(t)} ad_{b(t)}...$ are integrable . \qed

We can nw apply the procedure that we used for equations (\ref{lax2}): the exponential $exp_{Id_A + In_q(A)} $ exists and
$$ S_q (t) = Exp_{Id_A + In_q(A)}(ad_{P_q}) . S_q(0) .  \left(Exp_{Id_A + In_q(A)}(ad_{P_q})\right)^{-1}$$
is the unique solution to equation (\ref{symmetry2}) with initial value $S_q(0).$
\end{itemize}

\section{Appendix: a non regular Fr\"olicher Lie group}
Let $F$ be the vector space of smooth maps $C^\infty(]0;1[;\R).$ We equip $F$ with the following semi-norms:

 For each $(n,k) \in \N^* \times \N $, f $$||f||_{n,k} = \sup_{\frac{1}{n+1}\leq x \leq \frac{n}{n+1}} |f^{k}(x)|.$$

This is a Fr\'echet space.
Let 
$$\A = \{ f \in C^\infty(]0;1[;]0;1[)|\lim_{x \rightarrow 1}f(x)=1
 \wedge\lim_{x \rightarrow 0}f(x)=0 \}.$$
Finally, we set
 $$\D = \{ f \in \A | \inf_{x \in ]0;1[}f'(x) >0 \wedge \sup_{x \in ]0;1[}f'(x) >0\}.$$
 $\D$ is a set of diffeomorphisms of the open interval $]0;1[$ which is an (algebraic) group for composition of functions. Composition of maps is smooth for the topology on $\A$ induced by $C^\infty(]0;1[;\R)$. To see this, a well-know way to prove it is to pass to the diffeological setting (see e.g. \cite{Sou}
 for a review)
\begin{Proposition}
$\D$ is contractible.
\end{Proposition}
This proposition is trivial, but this will show us that our next results are not due to a pathology of the topology of $\D $ at the neighborhood of $\id .$
Unfortunately, $\D$ is not open in $\A.$
As a consequence, we are unable to prove that it is a Fr\'echet Lie group, sinse the actual knowledge on criteria to prove this is very restricted. However, considering the smooth diffeology induced on $\D$ by $\A,$ the inversion is smooth. As a consequence, $\D$ is a Fr\"olicher Lie group 
As a Fr\"olicher Lie group, it has a tangent space at $\id$ which is considered here as the Lie algebra.
Let $P\in \R[x]$ be a polynomial such that$$\left\{\begin{array}{l} P(0)=P(1)=0\\ \forall x\in ]0;1[, 0<P(x)<\min(x,1-x) \\ \sup_{x \in ]0;1[}|P'(x)|< 1 \end{array} \right. .$$
For example, one can take $$P(x)=\frac{x-x^2}{2}.$$ 
%$$C= 2\int_0^{+\infty} \frac{\ln(1+t^2)}{t^2} dt < +\infty,$$
We define 
$$\varphi(t,x)=\frac{P(x)t}{(1-P(x))t+P(x)},$$
$$c_t(x) = \left\{ \begin{array}{lr} x+\varphi(t,x) & \hbox{if } t\geq 0\\
 x - \varphi(-t,x) & \hbox{if } t<0. \end{array} \right.$$

Trivially, $\varphi$ is a smooth map defined on $\R_+ \times ]0;1[.$
So that, using the canonical isomorphism
$$ C^\infty(\R_+\times ]0;1[ ,\R) \sim C^\infty(\R_+,C^\infty(]0,1[,\R)),$$
we get that $$t\mapsto \varphi(t,.) \in C^\infty(\R_+,C^\infty(]0,1[,\R)),$$
and then $$ t \mapsto c_t\in C^1(\R,C^\infty(]0,1[,\R)).$$ 
We now want to show the following:
\begin{Theorem}\label{path}
The path $ t \mapsto c_t $ is of class $C^1(]-1;1[,\D).$ Moreover, 
$$ \partial_t c_t|_{t=0} = \mathbf{1}.$$
is a constant map.
\end{Theorem}

\begin{Lemma}
$$\forall t \in [0;1], c_t \in \D.$$
\end{Lemma}
\noindent
\textbf{Proof.} 
Let $x \in ]0;1[.$
The map $t \mapsto \varphi(t,x)$ is homographic. By standard results on homographic functions, it is strictly increasing on $[0;1].$ So that, if $t\in]-1;1[,$ $$-P(x) \leq \varphi(t,x) \leq P(x).$$  
As a consequence,
$$ \forall (t,x) \in ]-1;1[ \times ]0;1[, \quad 0 < x - {P(x)} < c_t(x) < x 
+ P(x)<1.$$
As a consequence, $l_0(c_t)=0$ and $l_1(c_t) = 1$ thus $c_t\in \A.$
Now, if $t \geq 0,$
\begin{eqnarray*}
 \partial_x \varphi(t,x) & = & \frac{tP'(x)(t+(1-t)P(x)) - tP(x)(1-t)P'(x) }{\left((1-P(x))t + P(x)\right)^2} \\
 & = &
  \frac{t^2P'(x)}{\left(t + P(x)(1-t)\right)^2}  .\end{eqnarray*}
We have $(1-t)P(x) \leq P(x) \leq 1,$ and then we get
$$\forall t \in ]-1;1[, \quad  1 - P'(x) \leq \partial_x c_t(x) \leq 1+P'(x).$$
Since $\sup_{x \in ]0;1[}|P'(x)|< 1,$ we get that $c_t$ is a diffeomorphism for each 
$t\in ]-1;1[$ which ends the proof.
  
\vskip 12pt
\noindent\textbf{Proof of Theorem \ref{path}.} Bb the previous Lemma, $t\mapsto c_t$ is a $C^1$ path with vales in $\D .$ Moreover, we only need to check that it is a $C^1-$ path for $t\geq 0.$ We have
\begin{eqnarray*}
 \partial_t \varphi(t,x) & = & \left(\frac{P(x)}{\left((1-P(x))t + P(x)\right)}\right)^2 \end{eqnarray*}
thus 
\begin{eqnarray*}
 \partial_t c_t(x) & = & \left(\frac{P(x)}{\left((1-P(x))t + P(x)\right)}\right)^2.
\end{eqnarray*} 
Thus $\partial_tc_t (x)|_{t=0} = 1.$
\vskip 12pt 
Now, we get the following:
\begin{Theorem} \label{nr}
There exists a smooth path $v $ on $T_{\id} \D$ such that no smooth path $g $ on $\D$ satisfies the equation $$\partial_t g \circ g^{-1} == v.$$
\end{Theorem}

\noindent
\textbf{Proof.}
Let $v$ be the constant path equal to $\mathbf{1}.$ Let $t \mapsto g_t$ be a solution of the last equation and let $x \in ]0;1[.$ Then we have, $\forall  g^{-1}(x) = y \in ]0;1[,$ $$\partial_t g_t (y) = 1,$$ so that the only possible solution is  the translation $g_t(x) = x+t.$ We have $$\forall t > 0, g_t \notin \A$$ so that $$\forall t>0, g_t \notin \D.$$

\end{document}